\documentclass[12pt]{article}%
\usepackage{amsfonts}
\usepackage{sw20bams}
\usepackage{amsmath}
\usepackage{amssymb}
\usepackage{graphicx}%
\providecommand{\U}[1]{\protect\rule{.1in}{.1in}}
\begin{document}

\title{Iterated Radical Expansions and Convergence}
\author{Steven Finch}
\date{November 5, 2024}
\maketitle

\begin{abstract}
We treat three recurrences involving square roots, the first of which arises
from an infinite simple radical expansion for the Golden mean, whose precise
convergence rate was made famous by Richard Bruce Paris in 1987. \ A
never-before-seen\ proof of an important formula is given. \ The other
recurrences are non-exponential yet equally interesting. \ Asymptotic series
developed for each of these two examples feature a constant, dependent on the
initial condition but otherwise intrinsic to the function at hand.

\end{abstract}

\footnotetext{Copyright \copyright \ 2024 by Steven R. Finch. All rights
reserved.}From the intricacies of nonlinear recurrences emerge a plethora of
constants. \ Our work builds on what we began in\ \cite{F1-radical,
F2-radical}.

\section{The Map $x\mapsto\sqrt{1+x}$}

As $k\rightarrow\infty$, the recurrence%
\[%
\begin{tabular}
[c]{llll}%
$x_{1}=1,$ &  & $x_{k}=\sqrt{1+x_{k-1}}$ & for $k\geq2$%
\end{tabular}
\ \ \ \
\]
approaches the Golden mean \cite{F3-radical}:%
\[%
\begin{array}
[c]{ccccc}%
\varphi=\dfrac{1+\sqrt{5}}{2}, &  & \varphi^{2}=\varphi+1, &  & \dfrac
{1}{\varphi}=\varphi-1
\end{array}
\]
from below and enjoys exponential convergence \cite{Pa-radical}:%
\[
0<\lim_{n\rightarrow\infty}\,(2\varphi)^{n}(\varphi-x_{n})=2\,%
{\displaystyle\prod\limits_{k=2}^{\infty}}
\,\frac{2\varphi}{\varphi+x_{k}}<\infty.
\]
We shall prove this formula using entirely elementary techniques. \ 

First, notice that $1\leq x_{k}<\varphi$ for all $k$ by induction ($x_{k}%
\geq1$ is obvious; supposing $1\leq x_{k-1}<\varphi$, we obtain%
\[
x_{k}=\sqrt{1+x_{k-1}}<\sqrt{1+\varphi}=\sqrt{1+\varphi^{2}-1}=\varphi
\]
). \ Now, writing $y_{k}=\varphi-x_{k}$, we have $y_{1}=\varphi-1=1/\varphi$,
$0<y_{k}\leq\varphi-1<\varphi^{2}$ and%
\begin{align*}
y_{k}  &  =\varphi-\sqrt{1+x_{k-1}}\\
&  =\varphi-\sqrt{1+\varphi-y_{k-1}}\\
&  =\varphi-\sqrt{\varphi^{2}-y_{k-1}}\\
&  <\dfrac{y_{k-1}}{\varphi}<\dfrac{y_{k-2}}{\varphi^{2}}<\dfrac{y_{k-3}%
}{\varphi^{3}}%
\end{align*}
because%
\[%
\begin{array}
[c]{ccccc}%
\dfrac{y_{k-1}}{\varphi}<\varphi, &  & \text{i.e.,} &  & \dfrac{y_{k-1}%
}{\varphi^{2}}<1
\end{array}
\]
hence%
\[%
\begin{array}
[c]{ccccc}%
0<\varphi-\dfrac{y_{k-1}}{\varphi}, &  & \text{i.e.,} &  & \dfrac{y_{k-1}^{2}%
}{\varphi^{2}}<y_{k-1}%
\end{array}
\]
hence%
\[
\varphi-\dfrac{y_{k-1}}{\varphi}=\sqrt{\varphi^{2}-2y_{k-1}+\dfrac{y_{k-1}%
^{2}}{\varphi^{2}}}<\sqrt{\varphi^{2}-y_{k-1}}%
\]
hence%
\[
\varphi-\sqrt{\varphi^{2}-y_{k-1}}<\dfrac{y_{k-1}}{\varphi};
\]
thus%
\[
y_{k}<\dfrac{y_{1}}{\varphi^{k-1}}=\frac{1}{\varphi^{k}}%
\]
for all $k$. \ 

From $x_{k-1}=\varphi-y_{k-1}$, we observe%
\[
x_{k}^{2}=1+x_{k-1}=1+\varphi-y_{k-1}=\varphi^{2}-y_{k-1}%
\]
hence%
\[
y_{k-1}=\varphi^{2}-x_{k}^{2}=(\varphi-x_{k})(\varphi+x_{k})=y_{k}%
(\varphi+x_{k})
\]
hence%
\[
\frac{2\varphi}{\varphi+x_{k}}=2\varphi\,\frac{y_{k}}{y_{k-1}}%
\]
hence%
\[%
{\displaystyle\prod\limits_{k=2}^{n}}
\,\frac{2\varphi}{\varphi+x_{k}}=%
{\displaystyle\prod\limits_{k=2}^{n}}
\,2\varphi\,\frac{y_{k}}{y_{k-1}}=(2\varphi)^{n-1\,}\frac{y_{n}}{y_{1}}%
=\frac{(2\varphi)^{n}}{2}y_{n}%
\]
hence%
\[
\frac{\varphi-x_{n}}{2}=\frac{y_{n}}{2}=\frac{1}{(2\varphi)^{n}}\,%
{\displaystyle\prod\limits_{k=2}^{n}}
\,\frac{2\varphi}{\varphi+x_{k}}%
\]
because $y_{1}=1/\varphi$; therefore%
\[
C=\lim_{n\rightarrow\infty}\,(2\varphi)^{n}\,\frac{\varphi-x_{n}}{2}=%
{\displaystyle\prod\limits_{k=2}^{\infty}}
\,\frac{2\varphi}{\varphi+x_{k}}%
\]
exists and is nonzero since%
\[%
{\displaystyle\sum\limits_{k=2}^{\infty}}
\left(  \frac{2\varphi}{\varphi+x_{k}}-1\right)  =%
{\displaystyle\sum\limits_{k=2}^{\infty}}
\,\frac{\varphi-x_{k}}{\varphi+x_{k}}<%
{\displaystyle\sum\limits_{k=2}^{\infty}}
\,y_{k}<%
{\displaystyle\sum\limits_{k=2}^{\infty}}
\,\frac{1}{\varphi^{k}}%
\]
converges. \ This completes the proof. \ The numerically-efficient expression
for%
\[
C=1.0986419643941564857346689...=\frac{1}{2}\left(
2.1972839287883129714693378...\right)
\]
as an infinite product did not appear in \cite{Pa-radical}, but was
subsequently found by Philippe Flajolet \& Paul Zimmermann several years later
\cite{Pf-radical}. \ No trace of their derivation has survived; the above
proof is new. \ Other calculations of $C$ are exhibited in \cite{Sn-radical,
S1-radical}.

In contrast, an analysis of the recurrence%
\[%
\begin{tabular}
[c]{llll}%
$x_{1}=0,$ &  & $x_{k}=\sqrt{\dfrac{1}{2}+\dfrac{1}{2}x_{k-1}}$ & for $k\geq2$%
\end{tabular}
\
\]
is trivial \cite{Au-radical}:%
\[%
\begin{array}
[c]{ccccc}%
x_{k}=\cos\left(  \dfrac{\pi}{2^{k}}\right)  , &  & x_{k}\rightarrow1^{-}, &
& \lim\limits_{k\rightarrow\infty}\,4^{k}(1-x_{k})=\dfrac{\pi^{2}}{2}.
\end{array}
\]
Finally, the recurrence%
\[%
\begin{tabular}
[c]{llll}%
$x_{1}=1,$ &  & $x_{k}=\sqrt{\,2+2x_{k-1}}$ & for $k\geq2$%
\end{tabular}
\
\]
approaches $1+\sqrt{3}$, but higher-order asymptotics (akin to
\cite{Pa-radical})\ await discovery.

\section{The Map $x\mapsto\dfrac{1+\sqrt{4x^{2}+1}}{2}$}

Given the recurrence%
\[%
\begin{array}
[c]{ccccc}%
x_{0}=1, &  & x_{k}=\dfrac{1+\sqrt{4x_{k-1}^{2}+1}}{2} &  & \text{for }k\geq1
\end{array}
\]
we start with asymptotics \cite{P1-radical}%
\[
x_{k}\sim\frac{1}{2}k+\frac{1}{4}\ln(k)+C
\]
valid as $k\rightarrow\infty$, for some constant $C=C(x_{0})$. \ On the basis
of numerical experimentation, we hypothesize that the next terms of the
asymptotic series must be of the form%
\begin{align*}
&  p_{1}\,\frac{\ln(k)}{k}+\frac{p_{0}}{k}+q_{2}\,\frac{\ln(k)^{2}}{k^{2}%
}+q_{1}\,\frac{\ln(k)}{k^{2}}+\frac{q_{0}}{k^{2}}+r_{3}\,\frac{\ln(k)^{3}%
}{k^{3}}+r_{2}\,\frac{\ln(k)^{2}}{k^{3}}\\
&  +r_{1}\,\frac{\ln(k)}{k^{3}}+\frac{r_{0}}{k^{3}}+s_{4}\,\frac{\ln(k)^{4}%
}{k^{4}}+s_{3}\,\frac{\ln(k)^{3}}{k^{4}}+s_{2}\,\frac{\ln(k)^{2}}{k^{4}}%
+s_{1}\,\frac{\ln(k)}{k^{4}}+\frac{s_{0}}{k^{4}}.
\end{align*}
The challenge is to express each coefficient $p_{j}$, $q_{j}$, $r_{j}$,
$s_{j}$ as a polynomial in $C$. \ To find these, we replace $k$ by $k+1$
everywhere:%
\[
\frac{1}{2}(k+1)+\frac{1}{4}\ln(k+1)+C+p_{1}\,\frac{\ln(k+1)}{k+1}%
+\cdots+\frac{s_{0}}{(k+1)^{4}}%
\]
and expand in powers of $k$ and $\ln(k)$:%

\[
\ln(k+1)\sim\frac{1}{k}-\frac{1}{2k^{2}}+\frac{1}{3k^{3}}-\frac{1}{4k^{4}%
}+\frac{1}{5k^{5}}-+\cdots,
\]%
\[
\frac{\ln(k+1)}{k+1}=\left(  \frac{1}{k}-\frac{1}{k^{2}}+\frac{1}{k^{3}}%
-\frac{1}{k^{4}}+\frac{1}{k^{5}}-+\cdots\right)  \ln(k)+\left(  \frac{1}%
{k^{2}}-\frac{3}{2k^{3}}+\frac{11}{6k^{4}}-\frac{25}{12k^{5}}+-\cdots\right)
,
\]%
\[
\frac{1}{k+1}\sim\frac{1}{k}-\frac{1}{k^{2}}+\frac{1}{k^{3}}-\frac{1}{k^{4}%
}+\frac{1}{k^{5}}-+\cdots,
\]%
\[
\frac{\ln(k+1)^{2}}{(k+1)^{2}}\sim\left(  \frac{1}{k^{2}}-\frac{2}{k^{3}%
}+\frac{3}{k^{4}}-\frac{4}{k^{5}}+\cdots\right)  \ln(k)^{2}+\left(  \frac
{2}{k^{3}}-\frac{5}{k^{4}}+\frac{26}{3k^{5}}-\cdots\right)  \ln(k)+\left(
\frac{1}{k^{4}}-\frac{3}{k^{5}}+\cdots\right)  ,
\]%
\[
\frac{\ln(k+1)}{(k+1)^{2}}\sim\left(  \frac{1}{k^{2}}-\frac{2}{k^{3}}+\frac
{3}{k^{4}}-\frac{4}{k^{5}}+-\cdots\right)  \ln(k)+\left(  \frac{1}{k^{3}%
}-\frac{5}{2k^{4}}+\frac{13}{3k^{5}}-+\cdots\right)  ,
\]%
\[
\frac{1}{(k+1)^{2}}\sim\frac{1}{k^{2}}-\frac{2}{k^{3}}+\frac{3}{k^{4}}%
-\frac{4}{k^{5}}+-\cdots,
\]%
\[
\frac{\ln(k+1)^{3}}{(k+1)^{3}}\sim\left(  \frac{1}{k^{3}}-\frac{3}{k^{4}%
}+\frac{6}{k^{5}}-+\cdots\right)  \ln(k)^{3}+\left(  \frac{3}{k^{4}}-\frac
{21}{2k^{5}}-+\cdots\right)  \ln(k)^{2}+\left(  \frac{3}{k^{5}}+-\cdots
\right)  \ln(k),
\]%
\[
\frac{\ln(k+1)^{2}}{(k+1)^{3}}\sim\left(  \frac{1}{k^{3}}-\frac{3}{k^{4}%
}+\frac{6}{k^{5}}-+\cdots\right)  \ln(k)^{2}+\left(  \frac{2}{k^{4}}-\frac
{7}{k^{5}}+-\cdots\right)  \ln(k)+\left(  \frac{1}{k^{5}}-+\cdots\right)  ,
\]%
\[
\frac{\ln(k+1)}{(k+1)^{3}}\sim\left(  \frac{1}{k^{3}}-\frac{3}{k^{4}}+\frac
{6}{k^{5}}-+\cdots\right)  \ln(k)+\left(  \frac{1}{k^{4}}-\frac{7}{2k^{5}%
}+-\cdots\right)  ,
\]%
\[
\frac{1}{(k+1)^{3}}\sim\frac{1}{k^{3}}-\frac{3}{k^{4}}+\frac{6}{k^{5}}%
-+\cdots,
\]%
\[
\frac{\ln(k+1)^{4}}{(k+1)^{4}}\sim\left(  \frac{1}{k^{4}}-\frac{4}{k^{5}%
}+-\cdots\right)  \ln(k)^{4}+\left(  \frac{4}{k^{5}}-+\cdots\right)
\ln(k)^{3},
\]%
\[
\frac{\ln(k+1)^{3}}{(k+1)^{4}}\sim\left(  \frac{1}{k^{4}}-\frac{4}{k^{5}%
}+-\cdots\right)  \ln(k)^{3}+\left(  \frac{3}{k^{5}}-+\cdots\right)
\ln(k)^{2},
\]%
\[
\frac{\ln(k+1)^{2}}{(k+1)^{4}}\sim\left(  \frac{1}{k^{4}}-\frac{4}{k^{5}%
}+-\cdots\right)  \ln(k)^{2}+\left(  \frac{2}{k^{5}}-+\cdots\right)  \ln(k),
\]%
\[
\frac{\ln(k+1)}{(k+1)^{4}}\sim\left(  \frac{1}{k^{4}}-\frac{4}{k^{5}}%
+-\cdots\right)  \ln(k)+\left(  \frac{1}{k^{5}}-+\cdots\right)  ,
\]%
\[
\frac{1}{(k+1)^{4}}\sim\frac{1}{k^{4}}-\frac{4}{k^{5}}+-\cdots.
\]
To avoid dealing with the radical in%
\[
2x_{k+1}-1=\sqrt{4x_{k}^{2}+1}%
\]
square both sides, obtaining%
\[
4x_{k+1}^{2}-4x_{k+1}+1=4x_{k}^{2}+1
\]
i.e.,%
\[
x_{k+1}^{2}-x_{k+1}=x_{k}^{2}.
\]
Upon rearrangement, relevant terms (in decreasing order of significance) of
$x_{k+1}^{2}-x_{k+1}$ become%

\begin{align*}
&  \left(  \frac{1}{8}-p_{1}+\frac{p_{0}}{2}+2p_{1}C+q_{1}\right)  \frac
{\ln(k)}{k}+\left(  -\frac{1}{8}+p_{1}-p_{0}+\frac{C}{2}+2p_{0}C+q_{0}\right)
\frac{1}{k}\\
&  +\left(  -\frac{p_{1}}{2}+p_{1}^{2}+r_{2}-2q_{2}+2q_{2}C+\frac{q_{1}}%
{2}\right)  \frac{\ln(k)^{2}}{k^{2}}\\
&  +\left(  -\frac{1}{16}+2p_{1}-\frac{p_{0}}{2}+2p_{1}p_{0}-2p_{1}%
C+r_{1}+2q_{2}-2q_{1}+2q_{1}C+\frac{q_{0}}{2}\right)  \frac{\ln(k)}{k^{2}}\\
&  +\left(  \frac{7}{48}-\frac{3p_{1}}{2}+\frac{3p_{0}}{2}+p_{0}^{2}-\frac
{C}{4}+2p_{1}C-2p_{0}C+r_{0}+q_{1}-2q_{0}+2q_{0}C\right)  \frac{1}{k^{2}}\\
&  +\left(  s_{3}-3r_{3}+2r_{3}C+\frac{r_{2}}{2}-q_{2}+2p_{1}q_{2}\right)
\frac{\ln(k)^{3}}{k^{3}}\\
&  +\left(  \frac{p_{1}}{2}-2p_{1}^{2}+s_{2}+3r_{3}-3r_{2}+2r_{2}C+\frac
{r_{1}}{2}+\frac{9q_{2}}{2}+2p_{0}q_{2}-4q_{2}C-q_{1}+2p_{1}q_{1}\right)
\frac{\ln(k)^{2}}{k^{3}}\\
&  +\left(  \frac{1}{24}-\frac{5p_{1}}{2}+2p_{1}^{2}+\frac{p_{0}}{2}%
-4p_{1}p_{0}+2p_{1}C+s_{1}+2r_{2}-3r_{1}+2r_{1}C+\frac{r_{0}}{2}-5q_{2}\right.
\\
&  \left.  \;\;\;\;\;+4q_{2}C+4q_{1}+2p_{0}q_{1}-4q_{1}C-q_{0}+2p_{1}%
q_{0}\right)  \frac{\ln(k)}{k^{3}}%
\end{align*}%
\begin{align*}
&  +\left(  -\frac{1}{8}+\frac{7p_{1}}{3}-\frac{7p_{0}}{4}+2p_{1}p_{0}%
-2p_{0}^{2}+\frac{C}{6}-3p_{1}C+2p_{0}C+s_{0}+r_{1}-3r_{0}\right. \\
&  \left.  \;\;\;\;\;+2r_{0}C+q_{2}-\frac{5q_{1}}{2}+2q_{1}C+\frac{7q_{0}}%
{2}+2p_{0}q_{0}-4q_{0}C\right)  \frac{1}{k^{3}}\\
&  +\left(  -4s_{4}+2s_{4}C+\frac{s_{3}}{2}-\frac{3r_{3}}{2}+2p_{1}r_{3}%
+q_{2}^{2}\right)  \frac{\ln(k)^{4}}{k^{4}}\\
&  +\left(  4s_{4}-4s_{3}+2s_{3}C+\frac{s_{2}}{2}+8r_{3}+2p_{0}r_{3}%
-6r_{3}C-\frac{3r_{2}}{2}+2p_{1}r_{2}+\frac{3q_{2}}{2}-6p_{1}q_{2}+2q_{2}%
q_{1}\right)  \frac{\ln(k)^{3}}{k^{4}}\\
&  +\left(  -\frac{p_{1}}{2}+3p_{1}^{2}+3s_{3}-4s_{2}+2s_{2}C+\frac{s_{1}}%
{2}-\frac{21r_{3}}{2}+6r_{3}C+\frac{15r_{2}}{2}+2p_{0}r_{2}-6r_{2}%
C-\frac{3r_{1}}{2}\right. \\
&  \left.  \;\;\;\;\;+2p_{1}r_{1}-\frac{31q_{2}}{4}+6p_{1}q_{2}-6p_{0}%
q_{2}+6q_{2}C+\frac{3q_{1}}{2}-6p_{1}q_{1}+q_{1}^{2}+2q_{2}q_{0}\right)
\frac{\ln(k)^{2}}{k^{4}}%
\end{align*}%
\begin{align*}
&  +\left(  -\frac{1}{32}+\frac{17p_{1}}{6}-5p_{1}^{2}-\frac{p_{0}}{2}%
+6p_{1}p_{0}-2p_{1}C+2s_{2}-4s_{1}+2s_{1}C+\frac{s_{0}}{2}+3r_{3}%
-7r_{2}+4r_{2}C\right. \\
&  \left.  \;\;\;\;\;+7r_{1}+2p_{0}r_{1}-6r_{1}C-\frac{3r_{0}}{2}+2p_{1}%
r_{0}+\frac{61q_{2}}{6}+4p_{0}q_{2}-10q_{2}C-\frac{13q_{1}}{2}\right. \\
&  \left.  \;\;\;\;\;+4p_{1}q_{1}-6p_{0}q_{1}+6q_{1}C+\frac{3q_{0}}{2}%
-6p_{1}q_{0}+2q_{1}q_{0}\right)  \frac{\ln(k)}{k^{4}}\\
&  +\left(  \frac{103}{960}-\frac{37p_{1}}{12}+p_{1}^{2}+\frac{23p_{0}}%
{12}-5p_{1}p_{0}+3p_{0}^{2}-\frac{C}{8}+\frac{11p_{1}C}{3}-2p_{0}%
C+s_{1}-4s_{0}+2s_{0}C\right. \\
&  \left.  \;\;\;\;\;+r_{2}-\frac{7r_{1}}{2}+2r_{1}C+\frac{13r_{0}}{2}%
+2p_{0}r_{0}-6r_{0}C-3q_{2}+2q_{2}C+\frac{29q_{1}}{6}+2p_{0}q_{1}\right. \\
&  \left.  \;\;\;\;\;-5q_{1}C-\frac{21q_{0}}{4}+2p_{1}q_{0}-6p_{0}q_{0}%
+6q_{0}C+q_{0}^{2}\right)  \frac{1}{k^{4}}.
\end{align*}
Performing an analogous substitution in $x_{k}$, the corresponding terms of
$x_{k}^{2}$ become%
\begin{align*}
&  \left(  \frac{p_{0}}{2}+2p_{1}C+q_{1}\right)  \frac{\ln(k)}{k}+\left(
2p_{0}C+q_{0}\right)  \frac{1}{k}+\left(  p_{1}^{2}+r_{2}+2q_{2}C+\frac{q_{1}%
}{2}\right)  \frac{\ln(k)^{2}}{k^{2}}\\
&  +\left(  2p_{1}p_{0}+r_{1}+2q_{1}C+\frac{q_{0}}{2}\right)  \frac{\ln
(k)}{k^{2}}+\left(  p_{0}^{2}+r_{0}+2q_{0}C\right)  \frac{1}{k^{2}}\\
&  +\left(  s_{3}+2r_{3}C+\frac{r_{2}}{2}+2p_{1}q_{2}\right)  \frac{\ln
(k)^{3}}{k^{3}}+\left(  s_{2}+2r_{2}C+\frac{r_{1}}{2}+2p_{0}q_{2}+2p_{1}%
q_{1}\right)  \frac{\ln(k)^{2}}{k^{3}}\\
&  +\left(  s_{1}+2r_{1}C+\frac{r_{0}}{2}+2p_{0}q_{1}+2p_{1}q_{0}\right)
\frac{\ln(k)}{k^{3}}+\left(  s_{0}+2r_{0}C+2p_{0}q_{0}\right)  \frac{1}{k^{3}}%
\end{align*}
\begin{align*}
&  +\left(  2s_{4}C+\frac{s_{3}}{2}+2p_{1}r_{3}+q_{2}^{2}\right)  \frac
{\ln(k)^{4}}{k^{4}}+\left(  2s_{3}C+\frac{s_{2}}{2}+2p_{0}r_{3}+2p_{1}%
r_{2}+2q_{2}q_{1}\right)  \frac{\ln(k)^{3}}{k^{4}}\\
&  +\left(  2s_{2}C+\frac{s_{1}}{2}+2p_{0}r_{2}+2p_{1}r_{1}+q_{1}^{2}%
+2q_{2}q_{0}\right)  \frac{\ln(k)^{2}}{k^{4}}\\
&  +\left(  2s_{1}C+\frac{s_{0}}{2}+2p_{0}r_{1}+2p_{1}r_{0}+2q_{1}%
q_{0}\right)  \frac{\ln(k)}{k^{4}}+\left(  2s_{0}C+2p_{0}r_{0}+q_{0}%
^{2}\right)  \frac{1}{k^{4}}.
\end{align*}
Matching coefficients, we obtain%
\[%
\begin{array}
[c]{ccc}%
p_{1}=\dfrac{1}{8}, &  & p_{0}=\dfrac{1}{2}C
\end{array}
\]
which are consistent with \cite{P2-radical} and%

\[%
\begin{array}
[c]{ccccc}%
q_{2}=-\dfrac{1}{32}, &  & q_{1}=-\left(  \dfrac{1}{4}C-\dfrac{1}{16}\right)
, &  & q_{0}=-\left(  \dfrac{1}{2}C^{2}-\dfrac{1}{4}C-\dfrac{1}{96}\right)  ,
\end{array}
\]
\[%
\begin{array}
[c]{ccccc}%
r_{3}=\dfrac{1}{96}, &  & r_{2}=\dfrac{1}{8}C-\dfrac{3}{64}, &  & r_{1}%
=\dfrac{1}{2}C^{2}-\dfrac{3}{8}C+\dfrac{1}{48},
\end{array}
\]%
\[
r_{0}=\frac{2}{3}C^{3}-\dfrac{3}{4}C^{2}+\dfrac{1}{12}C+\dfrac{7}{576},
\]%
\[%
\begin{array}
[c]{ccccc}%
s_{4}=-\dfrac{1}{256}, &  & s_{3}=-\left(  \dfrac{1}{16}C-\dfrac{11}%
{384}\right)  , &  & s_{2}=-\left(  \dfrac{3}{8}C^{2}-\dfrac{11}{32}%
C+\dfrac{5}{128}\right)  ,
\end{array}
\]%
\[%
\begin{array}
[c]{ccc}%
s_{1}=-\left(  C^{3}-\dfrac{11}{8}C^{2}+\dfrac{5}{16}C+\dfrac{1}{128}\right)
, &  & s_{0}=-\left(  C^{4}-\dfrac{11}{6}C^{3}+\dfrac{5}{8}C^{2}+\dfrac{1}%
{32}C-\dfrac{47}{5760}\right)
\end{array}
\]
which are new (as far as is known).

These fourteen parameter values allow us to estimate the constant $C$. \ Our
simple procedure involves computing $a_{10000000000}$ exactly via recursion,
setting this equal to our series (up to $s_{0}/k^{4}$) and then solving:%
\[
C=0.8232354508791921603541165....
\]
Note that the estimate $2\,C\approx1.6464707$ appears in \cite{S2-radical}.
\ We find the implicit representation%
\[
x_{k+1}(x_{k+1}-1)=x_{k}^{2}%
\]
to be intriguing: the left-hand side echoes the logistic map $\xi\,(1-\xi)$
but only somewhat:\ it is off by a sign. \ In the following section, a
comparable implicit representation leads to a surprising outcome.

\section{The Map $x\mapsto\dfrac{x+\sqrt{x^{2}+4}}{2}$}

Given the recurrence%
\[%
\begin{array}
[c]{ccccc}%
x_{0}=0, &  & x_{k}=\dfrac{x_{k-1}+\sqrt{x_{k-1}^{2}+4}}{2} &  & \text{for
}k\geq1
\end{array}
\]
we start with conjectured asymptotics (with no known theoretical basis)
\[
x_{k}\sim\sqrt{2k}-\frac{1}{4\sqrt{2}}\frac{\ln(k)}{k^{1/2}}-\frac{C}{k^{1/2}}%
\]
as $k\rightarrow\infty$, for some constant $C=C(x_{0})$. \ It is reasonable to
hypothesize that the next terms of this series are%
\[
p_{2}\,\frac{\ln(k)^{2}}{k^{3/2}}+p_{1}\,\frac{\ln(k)}{k^{3/2}}+\frac{p_{0}%
}{k^{3/2}}+q_{3}\,\frac{\ln(k)^{3}}{k^{5/2}}+q_{2}\,\frac{\ln(k)^{2}}{k^{5/2}%
}+q_{1}\,\frac{\ln(k)}{k^{5/2}}+\frac{q_{0}}{k^{5/2}}.
\]
To find $p_{j}$, $q_{j}$ we replace $k$ by $k+1$ everywhere:
\[
\sqrt{2(k+1)}-\frac{1}{4\sqrt{2}}\frac{\ln(k+1)}{(k+1)^{1/2}}+\frac
{C}{(k+1)^{1/2}}+p_{2}\,\frac{\ln(k+1)^{2}}{(k+1)^{3/2}}+\cdots+\frac{q_{0}%
}{(k+1)^{5/2}}%
\]
and expand in powers of $k$ and $\ln(k)$:%
\[
(k+1)^{1/2}\sim k^{1/2}+\frac{1}{2k^{1/2}}-\frac{1}{8k^{3/2}}+\frac
{1}{16k^{5/2}}-\frac{5}{128k^{7/2}}+\frac{7}{256k^{9/2}}-+\cdots,
\]%
\begin{align*}
\frac{\ln(k+1)}{(k+1)^{1/2}}  &  \sim\left(  \frac{1}{k^{1/2}}-\frac
{1}{2k^{3/2}}+\frac{3}{8k^{5/2}}-\frac{5}{16k^{7/2}}+\frac{35}{128k^{9/2}%
}-+\cdots\right)  \ln(k)\\
&  +\left(  \frac{1}{k^{3/2}}-\frac{1}{k^{5/2}}+\frac{23}{24k^{7/2}}-\frac
{11}{12k^{9/2}}+-\cdots\right)  ,
\end{align*}%
\[
\frac{1}{(k+1)^{1/2}}\sim\frac{1}{k^{1/2}}-\frac{1}{2k^{3/2}}+\frac
{3}{8k^{5/2}}-\frac{5}{16k^{7/2}}+\frac{35}{128k^{9/2}}-+\cdots,
\]%
\begin{align*}
\frac{\ln(k+1)^{2}}{(k+1)^{3/2}}  &  \sim\left(  \frac{1}{k^{3/2}}-\frac
{3}{2k^{5/2}}+\frac{15}{8k^{7/2}}-\frac{35}{16k^{9/2}}+-\cdots\right)
\ln(k)^{2}\\
&  +\left(  \frac{2}{k^{5/2}}-\frac{4}{k^{7/2}}+\frac{71}{12k^{9/2}}%
-+\cdots\right)  \ln(k)+\left(  \frac{1}{k^{7/2}}-\frac{5}{2k^{9/2}}%
+-\cdots\right)  ,
\end{align*}%
\[
\frac{\ln(k+1)}{(k+1)^{3/2}}\sim\left(  \frac{1}{k^{3/2}}-\frac{3}{2k^{5/2}%
}+\frac{15}{8k^{7/2}}-\frac{35}{16k^{9/2}}+\cdots\right)  \ln(k)+\left(
\frac{1}{k^{5/2}}-\frac{2}{k^{7/2}}+\frac{71}{24k^{9/2}}-\cdots\right)  ,
\]%
\[
\frac{1}{(k+1)^{3/2}}\sim\frac{1}{k^{3/2}}-\frac{3}{2k^{5/2}}+\frac
{15}{8k^{7/2}}-\frac{35}{16k^{9/2}}+-\cdots,
\]%
\begin{align*}
\frac{\ln(k+1)^{3}}{(k+1)^{5/2}}  &  \sim\left(  \frac{1}{k^{5/2}}-\frac
{5}{2k^{7/2}}+\frac{35}{8k^{9/2}}-+\cdots\right)  \ln(k)^{3}\\
&  +\left(  \frac{3}{k^{7/2}}-\frac{9}{k^{9/2}}+-\cdots\right)  \ln
(k)^{2}+\left(  \frac{3}{k^{9/2}}-+\cdots\right)  \ln(k),
\end{align*}%
\[
\frac{\ln(k+1)^{2}}{(k+1)^{5/2}}\sim\left(  \frac{1}{k^{5/2}}-\frac
{5}{2k^{7/2}}+\frac{35}{8k^{9/2}}-\cdots\right)  \ln(k)^{2}+\left(  \frac
{2}{k^{7/2}}-\frac{6}{k^{9/2}}+\cdots\right)  \ln(k)+\left(  \frac{1}{k^{9/2}%
}-\cdots\right)  ,
\]%
\[
\frac{\ln(k+1)}{(k+1)^{5/2}}\sim\left(  \frac{1}{k^{5/2}}-\frac{5}{2k^{7/2}%
}+\frac{35}{8k^{9/2}}-+\cdots\right)  \ln(k)+\left(  \frac{1}{k^{7/2}}%
-\frac{3}{k^{9/2}}+-\cdots\right)  ,
\]%
\[
\frac{1}{(k+1)^{5/2}}\sim\frac{1}{k^{5/2}}-\frac{5}{2k^{7/2}}+\frac
{35}{8k^{9/2}}-+\cdots.
\]
To avoid dealing with the radical in%
\[
2x_{k+1}-x_{k}=\sqrt{x_{k}^{2}+4}%
\]
square both sides, obtaining%
\[
4x_{k+1}^{2}-4x_{k+1}x_{k}+x_{k}^{2}=x_{k}^{2}+4
\]
i.e.,%
\[
x_{k+1}^{2}-1=x_{k+1}x_{k}%
\]
i.e.,%
\[
x_{k+1}-\frac{1}{x_{k+1}}=x_{k}.
\]
Upon rearrangement, the terms of $x_{k+1}-1/x_{k+1}$ involving either
$k^{5/2}$ or $k^{7/2}$ become most relevant. \ Performing an analogous
substitution in $x_{k}$, we match coefficients as before. \ The terms
containing $k^{5/2}$ give%
\[%
\begin{array}
[c]{ccccc}%
p_{2}=-\dfrac{1}{64\sqrt{2}}, &  & p_{1}=-\dfrac{4C-\sqrt{2}}{32}, &  &
p_{0}=-\dfrac{8\sqrt{2}C^{2}-8C+\sqrt{2}}{32}%
\end{array}
\]
and the terms containing $k^{7/2}$ give%
\[%
\begin{array}
[c]{ccccc}%
q_{3}=-\dfrac{1}{512\sqrt{2}}, &  & q_{2}=-\dfrac{3C-\sqrt{2}}{128}, &  &
q_{1}=-\dfrac{24\sqrt{2}C^{2}-32C+5\sqrt{2}}{256},
\end{array}
\]%
\[
q_{0}=-\frac{192C^{3}-192\sqrt{2}C^{2}+120C-11\sqrt{2}}{768}.
\]
It is astonishing that we have seen these coefficients before -- review our
analysis \cite{F2-radical} of the recurrence%
\[
\xi_{k+1}=\xi_{k}+\frac{1}{\xi_{k}}%
\]
which yields the \emph{identical} seven expressions -- although the signs
preceding each coefficient may differ. \ This is very surprising! \ While the
sequences behave distinctly, there is a hidden commonality in structure,
captured by the polynomials in $C$.

These seven parameter values allow us to estimate the constant $C$. \ Our
simple procedure involves computing $a_{10000000000}$ exactly via recursion,
setting this equal to our series (up to $q_{0}/k^{5/2}$) and then solving:%
\[%
\begin{array}
[c]{c}%
C=0.4117221539745403446660605....
\end{array}
\]
Note that the estimate $C/\sqrt{2}\approx0.291131527$ appears in
\cite{S3-radical}.

The seemingly arbitrary maps in Sections 2 \&\ 3 turn out to be curiously
linked: setting $y=1/x$, we have%
\[
\dfrac{1+\sqrt{4x^{2}+1}}{2}=x\,\dfrac{y+\sqrt{y^{2}+4}}{2}=x\,f(y),
\]%
\[
\,\dfrac{x+\sqrt{x^{2}+4}}{2}=x\,\dfrac{1+\sqrt{4y^{2}+1}}{2}=x\,g(y).
\]
Their appearance on a digital bulletin board \cite{S2-radical, S3-radical} by
different participants three years apart would suggest that there is no
connection. \ The mystics among us, however, might insist that there are no coincidences.

\section{Acknowledgements}

My proof in Section 1 is original, but it follows the outline of a related
proof that Robert Israel and Anthony Quas gave long ago (reproduced in
\cite{F1-radical}). I\ am thankful for helpful correspondence with Paul
Zimmermann and Michael Somos. My attempts to reach Dumitru Popa
\cite{P1-radical, P2-radical} have regrettably failed. \ In the statement of
Theorem 5 on page 21 of \cite{P2-radical}, the lead coefficient $1/2$ of the
$1/n$ term should be $1/4$. \ This particular theorem provides us rigorously
with the \textquotedblleft seed\textquotedblright\ (initial terms) from which
the asymptotic series in Section 2 grows. \ It does not cover the example in
Section 3 -- the reason is that the derivative $g^{\prime}(0)=0$ -- a
generalization over-and-beyond the archetypal scenario $f^{\prime}(0)\neq0$
will be needed to fill the gap.

\section{Addendum: $x\mapsto\sqrt{x(x+1)}$}

We've discovered that the expressions%
\[%
\begin{array}
[c]{ccccc}%
\dfrac{x+\sqrt{x^{2}+4}}{2} &  & \text{and} &  & x+\dfrac{1}{x}%
\end{array}
\]
enjoy a certain kindredship, in the sense that iterations based on these
functions possess strikingly similar asymptotic series. \ The same turns out
to be true for%
\[%
\begin{array}
[c]{ccccc}%
\dfrac{1+\sqrt{4x^{2}+1}}{2} &  & \text{and} &  & \sqrt{x(x+1)}.
\end{array}
\]
Consider the recurrence $x_{0}=1$, $x_{k}=\sqrt{x_{k-1}(x_{k-1}+1)}$ for
$k\geq1$. \ We find%
\begin{align*}
x_{k}  &  \sim\frac{1}{2}k-\frac{1}{4}\ln(k)-C+\frac{1}{8}\frac{\ln(k)}%
{k}+\frac{C}{2}\frac{1}{k}+\frac{1}{32}\frac{\ln(k)^{2}}{k^{2}}+\left(
\frac{1}{4}C-\frac{1}{16}\right)  \frac{\ln(k)}{k^{2}}\\
&  +\left(  \dfrac{1}{2}C^{2}-\dfrac{1}{4}C-\dfrac{1}{96}\right)
\frac{1^{\backslash}}{k^{2}}+\frac{1}{96}\frac{\ln(k)^{3}}{k^{3}}+\left(
\frac{1}{8}C-\frac{3}{64}\right)  \frac{\ln(k)^{2}}{k^{3}}\\
&  +\left(  \dfrac{1}{2}C^{2}-\dfrac{3}{8}C+\dfrac{1}{48}\right)  \frac
{\ln(k)}{k^{3}}+\left(  \frac{2}{3}C^{3}-\dfrac{3}{4}C^{2}+\dfrac{1}%
{12}C+\dfrac{7}{576}\right)  \frac{1}{k^{3}}\\
&  +\dfrac{1}{256}\frac{\ln(k)^{4}}{k^{4}}+\left(  \dfrac{1}{16}C-\dfrac
{11}{384}\right)  \frac{\ln(k)^{3}}{k^{4}}+\left(  \dfrac{3}{8}C^{2}%
-\dfrac{11}{32}C+\dfrac{5}{128}\right)  \frac{\ln(k)^{2}}{k^{4}}\\
&  +\left(  C^{3}-\dfrac{11}{8}C^{2}+\dfrac{5}{16}C+\dfrac{1}{128}\right)
\frac{\ln(k)}{k^{4}}+\left(  C^{4}-\dfrac{11}{6}C^{3}+\dfrac{5}{8}C^{2}%
+\dfrac{1}{32}C-\dfrac{47}{5760}\right)  \frac{1}{k^{4}}%
\end{align*}
and, unlike the alternating blocks of positive \& negative signs in Section 2,
here all signs are positive (except for the earliest $\ln(k)$ \&\ $C$ terms).
\ More on this phenomenon is given in \cite{F4-radical}. \ Using the same
procedure as before, we obtain%
\[
C=-1.1751774424585571398132856....
\]
The commonality in structure between two functions does not assist in the
numerical calculation of constants. \ No relationship between $0.823...$ and
$-1.175...$ is observed (nor was any expected).

\end{document}